\documentclass[a4paper,12pt,reqno]{amsart}

\usepackage[all]{xy}
\usepackage[english]{babel}
\usepackage{amssymb}
\usepackage{graphicx}
\usepackage{enumerate}
\usepackage{amsmath,amsthm,amssymb,xypic}
\usepackage{xcolor}
\usepackage{palatino}
\usepackage{geometry}
\newcommand{\p}{\mathbb{P}}
\newcommand{\pp}{\mathbb{P}}

\newcommand{\N}{\mathbb{N}}

\newcommand{\C}{\mathbb{C}}
\newcommand{\vu}{\varnothing}
\newcommand{\I}{\mathcal{I}}

\newcommand{\G}{\mathbb{G}}
\newcommand{\T}{\mathbb{T}}
\newcommand{\f}{\varphi}
\newcommand{\LL}{\mathcal{L}}
\newcommand{\oo}{\mathcal{O}}
\newcommand{\mul}{\mathop{\rm mult}\nolimits}
\newcommand{\Bs}{\mathop{\rm Bs}\nolimits}

\newcommand{\SV}{\mathop{\rm SV}\nolimits}
\newcommand{\codim}{\mathop{\rm codim}\nolimits}

\newcommand{\V}{\mathop{\rm V}\nolimits}
\newcommand{\Sec}{\mathop{\rm Sec}\nolimits}
\newcommand{\SSec}{\mathop{\rm \mathbb{S}ec}\nolimits}

\theoremstyle{plain}

\theoremstyle{definition}
\newtheorem{thm}{Theorem}
\newtheorem{pro}[thm]{Proposition}
\newtheorem{lem}[thm]{Lemma}

\newtheorem{dfn}[thm]{Definition}
\newtheorem{rmk}[thm]{Remark}

\newcommand{\Fra}[1]{{
#1}}

\begin{document}

\title[Generic identifiability of pairs of ternary forms]{Generic identifiability of pairs of ternary forms}

\author[V. Beorchia]{Valentina Beorchia}
\address{Department of Mathematics and Geosciences\\ University of Trieste\\ via Valerio 12/1 \\ 34127 Trieste, Italy,
ORCID 0000-0003-3681-9045.}
\email{beorchia@units.it}

\author[F. Galuppi]{Francesco Galuppi}
\address{Institute of Mathematics of the Polish Academy of Sciences\\ ulica Śniadeckich 8\\00-656 Warszawa, Poland, ORCID 0000-0001-5630-5389.}

\email{fgaluppi@impan.pl}

 

\thanks{The authors 
acknowledge the MIUR Excellence Department Project awarded to the Department of Mathematics and Geosciences, University of Trieste. Beorchia is member of GNSAGA of INdAM and is  supported  by  italian MIUR funds, PRIN  project {\it Moduli  Theory and Birational Classification} (2017), P.i. U.  Bruzzo. Galuppi is supported by the National Science Center, Poland, project ``Complex contact manifolds and geometry of secants'', 2017/26/E/ST1/00231. The authors wish to thank the anonymous referee for useful comments on the first version of this paper.}

\begin{abstract} 
We prove that two general ternary forms \Fra{of degrees $c\le d$} are simultaneously identifiable only in the classical cases $(c,d)=(2,2)$ and $(c,d)=(2,3)$.
We translate the problem into the study of a certain linear system on a projective bundle on the plane, and we apply techniques from projective and birational geometry to prove that the associated map is not birational. \end{abstract}
\maketitle

\section{Introduction}\label{sec: intro}Many questions in mathematics are called \textquotedblleft Waring problem\textquotedblright, after the name of the 18th century English mathematician Edward Waring. He was interested in decomposing a natural number as a sum of powers. Since then, Waring problems have to do with additive decompositions of different mathematical objects. For instance, given a degree $d$ form, or equivalently, a homogeneous polynomial $f$, we can decompose
\[f=\ell_1^d+\ldots+\ell_k^d\]
as a sum of powers of linear forms. In a similar way, one can decompose a tensor as a sum of rank one tensors. Waring decompositions raise a huge interest in many different areas, as well illustrated in the textbook \cite[Section 1.3]{Lan}.

There are several different questions that we may ask. For instance, what is the smallest possible number of summands, called the \emph{Waring rank} of $f$, or how to compute all the 
 decompositions. When they are infinitely many, one can study the variety parametrizing them. When they are finitely many, one may ask to bound their number. All these questions are widely open in their generality. One interesting problem is to understand when the decomposition is unique. In this case, we say that $f$ is  $k$-\emph{identifiable}.  Identifiability is a desirable property with many applications, as it gives a canonical form for $f$. Examples range from Signal processing to Complexity theory, from Philogenetics to Algebraic statistics. A complete list of applications would be far too long, so we refer the reader to \cite[Section 1.3]{Lan}.

\Fra{When we work over the complex field, there is a dense open subset of the space of polynomials where all elements have the same rank - called the \textit{generic rank} - and the same number of decompositions. In this paper we will use the words "generic" or "general" for properties that hold almost everywhere - more precisely, on a dense open subset.} Classically, the problem was to classify all pairs $(n,d)$ such that the general $f\in\C[x_0,\dots,x_n]_d$ is identifiable. General identifiability is expected to be a rare phenomenon. A few classical cases were known to be generically identifiable since the work of Hilbert and Sylvester. It took more than a century to get new results in this direction \cite{M1,M2}, and the full classification has been completed in \cite{GM}. It turns out that there are infinitely many generically identifiable cases for binary forms, while there are only two sporadic cases for polynomials in three or more variables.

In this paper we focus on the version of the Waring problem concerning pairs of polynomials. It is a classical result that two general quadratic forms $f,g\in\C[x_0,\dots,x_n]_2$ can be simultaneously diagonalized. More precisely, there exist linear forms $\ell_1,\ldots, \ell_{n+1}$ and scalars $\lambda_1,\ldots, \lambda_{n+1}$ such that
\begin{equation}\label{eq:simultaneous diagonalization}
\begin{cases}
f=\ell_1^2+\ldots+\ell_{n+1}^2\\
g=\lambda_1\ell_1^2+\ldots+\lambda_{n+1}\ell_{n+1}^2.
\end{cases}
\end{equation}
A canonical form \eqref{eq:simultaneous diagonalization} with $n+1$ summands 
is unique if and only if the discriminant of the pencil $\langle f,g\rangle$ does not vanish, hence the general pair of quadratic forms has only one simultaneous diagonalization.

We generalize decomposition \eqref{eq:simultaneous diagonalization} to pairs of forms of any degrees. For symmetry reasons, it is convenient not to distinguish $f$ from $g$, so we will allow coefficients in the decomposition of $f$ as well. 

\begin{dfn}\label{simdec} Let $f\in\C[x_0,\ldots,x_n]_c$ and $g\in\C[x_0,\ldots,x_n]_d$ be two homogeneous polynomials. A \emph{Waring decomposition} of  $(f,g)$ consists of linear forms $ \ell_1, \ldots, \ell_k \in \C[x_0,\ldots,x_n]_1 $ and scalars $ \lambda_1,\ldots,\lambda_k,\mu_1,\ldots,\mu_k \in \mathbb{C}$ such that
	\begin{equation}\label{eq:simultaneous decomposition}
	\begin{cases}
	f=\lambda_1\ell_1^c+\ldots+\lambda_k\ell_k^c\\
	g=\mu_1\ell_1^d+\ldots+\mu_k\ell_k^d.
	\end{cases}
	\end{equation}
This kind of decompositions are also called \emph{simultaneous decompositions}. Due to the presence of the scalars $ \lambda_1,\ldots,\lambda_k$ and $\mu_1,\ldots,\mu_k$, each linear form depends essentially only on $n$ conditions, so the we can regard \eqref{eq:simultaneous decomposition} as a polynomial system with $\binom{c+n}{n}+\binom{d+n}{n}$ equations - given by the data $f$ and $g$ - and $k(n+2)$ unknowns - namely, the scalars $ \lambda_1,\ldots,\lambda_k,\mu_1,\ldots,\mu_k$ and the linear forms $ \ell_1, \ldots, \ell_k$. 
We consider two decompositions of $(f,g)$ to be equal if they differ just by the order of the $k$ summands. The \emph{Waring rank}, or \emph{simultaneous rank} of	$(f,g)$ is the minimum number $k$ such that there exists a simultaneous decomposition \eqref{eq:simultaneous decomposition} with $k$ summands.
A pair is \emph{$k$-simultaneously identifiable}, or simply \emph{$k$-identifiable} if it admits a unique simultaneous decomposition with $k$ summands.
\end{dfn}
This problem fascinates mathematicians since a long time ago. In \cite{London}, London proved that the rank of two general ternary cubics is 6, instead of the expected number 5. Later, Scorza described London's result from a different perspective in \cite{Scorza}\Fra{ - see \cite[Theorem 4.1]{CarChi03} for a modern reference}. In \cite{terracinicoppiediternarie}, Terracini computed the simultaneous rank of two general ternary forms of the same degree. As described in \cite[Section 2.2]{AGMO}, we can rephrase these results in modern language: if we call $\SV_{1\times n}^{1,d}$ the Segre-Veronese embedding of $\p^1\times\p^n$ via the complete linear system of divisors of bidegree $(1,d)$, then London proved that $\SV_{1\times 2}^{1,3}$ is 5-defective, while Terracini showed that $\SV_{1\times 2}^{1,d}$ is not defective for $d\neq 3$. We had to wait until the work \cite{BalBerCat12} to have a classification of all defective $\SV_{1\times n}^{1,d}$, so now we know the rank of two general forms of the same degree in any number of variables. However, when the two degrees are different, the problem of computing the rank is not yet solved. 

Papers like \cite{LRA}, \cite{REF} and  \cite{MOUR} use simultaneous diagonalization of 
symmetric matrices 
to bound the rank of a third-order symmetric tensor
. We believe that simultaneous decompositions like \eqref{eq:simultaneous decomposition} 
could be applied in a similar way to study the rank of higher-order symmetric tensors. 

In this paper we focus on identifiability. The guiding problem is the classification of all triples $(n,c,d)$ such that the general pair of forms of degrees $c$ and $d$ in $n+1$ variables is identifiable. In \cite[Section 5]{CR}, Ciliberto and Russo solve the case $n=1$ of binary forms. They work in a slightly different language and phrase their statement in terms of geometric properties of a rational normal scroll. Their result applies to tuples of binary forms, not just pairs - see also \cite[Theorem 1.3]{AGMO}. Roughly speaking, the general pair of binary forms is identifiable, as long as the polynomial system \eqref{eq:simultaneous decomposition} is square and $d$ is not too large compared to $c$. This reminds us what happens for the Waring problem for one polynomial: while generic identifiability is rare, there are infinitely many cases for binary forms.
	
The situation changes when there are more than two variables. As widely expected, for $n\ge 2$ generic identifiability is very uncommon. In more than a century, mathematicians have found only two instances in which the general pair of forms is identifiable. The case \eqref{eq:simultaneous diagonalization} of two general quadrics goes back at least to Weierstrass and it is generically identifiable in any number of variables. Beside that, there is the case of a general plane conic and a cubic, studied by Roberts 
	in \cite{Roberts} \Fra{and revisited in \cite[Theorem 10.2]{OS}}. The challenge is to prove that those are the only cases, or to find new exceptional ones. In this paper we solve the problem for ternary forms.
	
	\begin{thm}
		\label{thm: main} Let $c$ and $d$ be positive integers such that $c\le d$.	The general pair of ternary forms of degrees $c$ and $d$ is identifiable if and only if $(c,d)\in\{(2,2), (2,3)\}$.
	\end{thm}
In the special case $d=c+1$, Theorem \ref{thm: main} has been proven in \cite[Theorem 5.1]{AGMO}. Despite the algebraic statement, our approach relies on projective and birational geometry. In Section \ref{section: problem reduction}, we underline the tight connection between decompositions and secant varieties and we translate the problem into a question about the degree of a certain rational map. Namely, the set of decompositions of a pair $(f,g)$ 
is the fiber of the 
secant map of the projective bundle $X=\p(\oo_{\p^2}(c)\oplus\oo_{\p^2}(d))$. In order to disprove identifiability, we show that the map is not birational. As in \cite[Section 5]{AGMO}, the first step is to bound the degree of such a map with the degree of a certain linear projection of $X$. Then we degenerate the associated linear system and we restrict it to a suitable subvariety to prove that such a map cannot be birational.

When performing this kind of degenerations, it is common to bump into some artithmetic obstructions. We overcome this obstacle by distinguishing two cases and give two different arguments, in Sections \ref{section: first case} and \ref{sec:secondo caso} respectively. For this reason the proof of Theorem \ref{thm: main} is split in two parts, namely Propositions \ref{pro: dimostrazione nel primo caso} and \ref{pro: dimostrazione nel secondo caso}. 

\section{Geometric setup}\label{section: preliminaries}

\Fra{Secant varieties are a classical construction that dates back to the Italian school of algebraic geometry at the end of 19th century. In this section we recall the definition of the secants of a variety $V$ embedded in some projective space $\p^N$. We will use such a topic in the case in which $V$ is the variety parametrizing pairs of polynomials of simultaneous rank 1.

We work over the complex field $\C$. }
Let $\G(k-1,N)$ be the Grassmannian of $(k-1)$-linear spaces in $\p^N$.
Let $V\subset\p^{N}$ be a nondegenerate irreducible variety of dimension $ n $ and let
\[\Gamma_{k}(V)=\overline{\{(x_1,\ldots,x_k,L)\in V\times\ldots\times V\times\G(k-1,N)\mid L=\langle x_1,\ldots,x_k\rangle\}}.
\]
Observe that $\Gamma_{k}(V)$ is birational to $V\times\ldots\times V$, therefore it is irreducible of dimension $kn$. Let $\pi_2:\Gamma_{k}(V)\to\G(k-1,N)$ be the projection onto the last factor and set
\[S_k(V)=\pi_2(\Gamma_{k}(V))=\overline{\{L\in \G(k-1,N)\mid L\mbox{ is spanned by $k$ points of }V \}}.\]
\Fra{Thanks to the Trisecant lemma - see for instance \cite[Proposition 2.6]{trisecant} - t}he general $L\in S_k(V)$ meets $V$ in exactly $k$ points, so the general fiber of $\pi_2$ has dimension zero. Hence $$\dim S_k(V)=\dim(\Gamma_{k}(V))=kn.$$ We are ready to define the secant varieties of $V$.

\begin{dfn}\label{def:secant}
Let $V\subset\p^{N}$ be a nondegenerate irreducible variety. The {\it abstract $k$-secant variety} of $V$ is
$$\Sec_k(V)=\overline{\{
(x,L)\in \p^N\times\G(k-1,N)\mid x\in L\mbox{ and }L\in S_k(V)\}}.$$
Let $p_1:\Sec_k(V)\to\p^N$ and $p_2:\Sec_k(V)\to S_k(V)$ be the two projections. The general fiber of $p_2$ is a linear space of dimension $k-1$, therefore $$\dim\Sec_k(V)=\dim S_k(V)+k-1=kn+k-1.$$
The {\it $k$-secant variety} of $V$ is
$$\SSec_k(V)=p_1(\Sec_k(V))=\overline{\bigcup_{x_1,\dots,x_k\in V}\langle x_1,\dots,x_k\rangle}\subset\p^N.$$
By construction, $\dim\SSec_k(V)\le\min\{\dim\Sec_k(V),N\}=\min\{kn+k-1,N\}$. The variety $V$ is called
$k$-\emph{defective} if $$
\dim\SSec_k(V)< \min\{kn+k-1,N\}.
$$ 
\end{dfn}

A classical result about secant varieties is Terracini's lemma, proven in \cite{terracini} - see \cite[Section 2]{dale} for a modern reference.
\begin{lem}[Terracini]\label{lem: Terracini} Let $V\subset\p^{N}$ be a nondegenerate irreducible variety. 
	If $x_1,\ldots, x_r\in V$ are in general position and $z\in\langle x_1,\ldots, x_{r}\rangle$
	is a general point, then the embedded tangent space to $\SSec_r(V)$ at $z$ is
	$$\T_z\SSec_r(V) = \langle \T_{x_1}V,\ldots, \T_{x_{r}}V\rangle.$$
\end{lem}

Terracini's lemma allows us to link the study of secant varieties of $V$ to the study of linear systems of divisors of $V$ with imposed singularities. \Fra{Indeed, 
let $\LL$ be the linear system on $V$ of hyperplane sections
. Then
\begin{align*}
\codim\SSec_r(V) & =\codim \T_z\SSec_r(V)\\
&=\dim\{H\subset\p^N\mid H\mbox{ is a hyperplane and } H\supset\T_z\SSec_r(V)\}\\
&=\dim\{H\cap V\mid H\mbox{ is a hyperplane and } H\supset\T_z\SSec_r(V)\}\\
&=\dim\{H\cap V\mid H\mbox{ is a hyperplane and } H\supset\T_{x_i}V\mbox{ for every } i\in\{1,\dots,r\}\}\\
&=\dim\{D\in\LL\mid D\mbox{ is singular at }x_i\mbox{ for every } i\in\{1,\dots,r\}\}.
\end{align*}
}
Let us recall some definitions we will use when dealing with linear systems.

\begin{dfn}
Let $p$ be a point on a variety $V$ defined by the ideal $I_p$ and let $m\in\N$. The \emph{point of multiplicity $m$} supported at $p$ is the $0$-dimensional subscheme of $V$ defined by $I_p^m$. A point of multiplicity 2 is also called a \emph{double point}.
\end{dfn}

Since 
we are going to work with systems of plane curves, we introduce the notation we use. \Fra{We denote by
\[L_2^d(2^r,1^t)\]
the vector space of degree $d$ ternary forms vanishing at $r$ \Fra{fixed} general points with multiplicity at least 2 and vanishing at $t$ further \Fra{points in general position}. Its projectivization
\[\LL_2^d(2^r,1^t)=\p(L_2^d(2^r,1^t))\]
is the linear system of degree $d$ plane curves singular at $r$ fixed  general points and containing $t$ base points in general position. \Fra{More generally, given a linear system $\LL=\p(L)$ on a variety $V$, we denote by $\LL(2^r,1^t)$ the linear subsystem of $\LL$ consisting of divisors singular at $r$ fixed general points and vanishing at $t$ further points in general position. When $\LL$ is the linear system embedding $V\subset\p^N$, that is, the linear system of hyperplane sections of $V$, Lemma \ref{lem: Terracini} tells us that
\begin{equation}\label{eq: sistemi lineari vs secanti}
    \codim\SSec_r(V)=\dim L(2^r).
\end{equation}
It is natural to ask for the dimension of the vector space $L_2^d(2^r, 1^t)$. When the base points have arbitrary multiplicity, that can be very hard to compute. However, when the multiplicities are at most 2, the answer is provided by the celebrated Alexander-Hirschowitz' theorem
. Such result is a landmark in the fields and holds in a projective space of any dimension, but here we only recall the weaker version on the plane, that is known at least since \cite{Campbell}.
\begin{thm}\label{thm: AH}
The dimension of $L_2^d(2^r)$ is $\max\left\lbrace 0,\binom{d+2}{2}-3r\right\rbrace$, unless $(d,r)\in\{(2,2),(4,5)\}$. In these exceptional cases we have $\dim L_2^2(2^2)=\dim L_2^4(2^5)=1$.
\end{thm}
Since simple points in general position always impose independent conditions, Theorem \ref{thm: AH} gives us a formula to compute $\dim L_2^d(2^r,1^t)$.
}} When handling such systems, we will need some control on their  singularities. The following is a slight generalization of a result proven in \cite{AC}.

\begin{lem}
\label{lemma: arbarellocornalba} Let $d$, $r$ and $t$ be natural numbers such that $(d,r,t)\neq (6,9,0)$. If
\[
\binom{d-1}{2}\ge r\mbox{ and } \binom{d+2}{2}-1\ge 3r+t,
\]
then the general element of $\LL_2^d(2^r,1^t)$ is irreducible, has exactly $r$ ordinary double points and it is smooth elsewhere.
\begin{proof}
\Fra{By Theorem \ref{thm: AH}, t}he hypothesis $\binom{d+2}{2}-1\ge 3r+t$ guarantees that $\LL_2^d(2^r,1^t)$ is not empty. 
We argue by induction on $t$. For $t=0$, the claim is proven in \cite[Theorem 3.2]{AC}. Now assume that $t\ge 1$. \Fra{The system $\LL_2^6(2^9,1^t)$ is empty by Theorem 6, so }we assume that $(d,r)\neq (6,9)$. By induction hypothesis, $\LL_2^d(2^r,1^{t-1})$ contains a nonempty open subset $U$ consisting of irreducible curves with exactly $r$ ordinary double points. Imposing a further simple point in the base locus corresponds to taking a hyperplane section of $\LL_2^d(2^r,1^{t-1})\subset \p (\C [x_0,x_1,x_2]_d)$; if the hyperplane is general, the intersection with $U$ is nonempty.
%
\end{proof}
\end{lem}

There is one last thing we need to recall before we move to the next section. Let $V$ be a projective variety and let $Z$ be a subvariety of $V$. Consider a linear system $\LL=\p(L)$ on $V$
and let $\LL_{Z}=\p(L_Z)$ be the complete linear system on $Z$ given by $\p(H^0\oo_Z(D_{|Z}))$ for a general element $D\in\LL$.
Then there is an exact sequence of vector spaces
\begin{equation}\label{exact_sequence}
0 \rightarrow L\cap I_{Z} \rightarrow L \rightarrow L_Z.
\end{equation}
The image of the rightmost map is denoted by $L_{|Z}$ and its projectivization is $\LL_{|Z}=\{D_{|Z}\mid D\in\LL\}$.
Sequence \eqref{exact_sequence} is called \textit{Castelnuovo exact sequence}.

\section{Problem reduction}\label{section: problem reduction}
In this section we formalize the problem and we present some simplifications. Let us start by considering the Waring problem for one polynomial. Degree $d$ ternary forms of Waring rank 1 are parametrized by the \emph{$d$-Veronese surface}
\[\V_2^d=\{[\ell^d]\mid\ell\in\C[x_0,x_1,x_2]_1\}\subset\p^{\binom{n+2}{2}-1}=\p(\C[x_0,x_1,x_2]_d)\]
embedded in the space of all forms of degree $d$. Then the rank of the general ternary form is\[\min\{r\in\N\mid\SSec_r(\V_2^d)=\p^{\binom{d+2}{2}-1}\}.\]
Thanks to Lemma \ref{lem: Terracini} and equation \eqref{eq: sistemi lineari vs secanti}, the latter equals
\[\min\{r\in\N\mid\dim L_2^d(2^r)=0\}.
\]
This allows us to use geometric techniques to study Waring decompositions. For instance, we can apply Theorem \ref{thm: AH} to compute the rank of the general ternary form of degree $d$. Moreover, the set of decompositions with $k$ summands of a polynomial $f\in\p^{\binom{d+2}{2}-1}$ is the fiber of the secant map $p_1:\Sec_k(\V_2^d)\to\p^{\binom{d+2}{2}-1}$ over $f$ - see Definition \ref{def:secant}.

We consider a similar construction for simultaneous decompositions. Let $c$ and $d$ be positive integers such that $c\le d$. The variety parametrizing pairs of polynomials of degrees $c$ and $d$ and simultaneous rank $1$ is 
$$X=\{[a_1\ell^c,a_2\ell^d]\mid\ell\in\C[x_0,x_1,x_2]_1\mbox{ and }a_1,a_2\in\C\}\subset\p(\C[x_0,x_1,x_2]_c\oplus\C[x_0,x_1,x_2]_d).$$
The map $$\pi:X\to\p^2=\p(\C[x_0,x_1,x_2]_1)$$ sending $[a_1\ell^c,a_2\ell^d]$ to $[\ell]$
gives $X$ the structure of a projective bundle over the plane, where each fiber is isomorphic to $\p^1$. Therefore $X$ is a threefold. For the basic definitions and notions concerning projective bundles, we refer the reader to \cite [Chapter II.7]{hartshorne} or \cite[Chapter 9]{3264}. It turns out that $X$ is the projectivization of the following rank two vector bundle on $\p^2$:
\[X=\p(\oo_{\p^2}(c)\oplus\oo_{\p^2}(d))\subset\p^N=\p(\C[x_0,x_1,x_2]_c\oplus\C[x_0,x_1,x_2]_d),
\]
where $N=\binom{c+2}{2}+\binom{d+2}{2}-1$. The immersion $X\subset\p^N$ is the tautological embedding. Namely, the very ample divisor $T_X$ associated with the immersion - called the \textit{tautological divisor} - is the only divisor class on $X$ such that $$\pi_*\oo_X(T_X)=\oo_{\p^2}(c)\oplus\oo_{\p^2}(d).$$
We recall that $T_X$ is \textit{unisecant}, that is, each fiber of $\pi$ intersects $T_X$ in exactly one point.


The situation is similar to the case of the Veronese variety. According to \cite[Section 2.2]{AGMO}, the set of decompositions with $k$ summands of $(f,g)\in\p^N$ is the fiber of the secant map $p_1:\Sec_k(X)\to \p^N$ over $(f,g)$. For the generic pair of polynomials, in order to have finitely many decompositions, we have to assume that $\dim\Sec_k(X)=\dim \p^N$. This leads to the following definition.

\begin{dfn} Let $c$ and $d$ be positive integers such that $c\le d$.\color{black} \label{def: definizione di k}
We say that $(c,d)$ is a \emph{perfect case} if  there exists $k\in\N$ such that $\dim\Sec_k(X)=\dim \p^N$. This is equivalent to $3k+k-1=\binom{c+2}{2}+\binom{d+2}{2}-1$, so $(c,d)$ is perfect if and only if $\binom{c+2}{2}+\binom{d+2}{2}$ is a multiple of 4. In this case 
\begin{equation}
    \label{equation:definition of k}k=\frac{\binom{c+2}{2}+\binom{d+2}{2}}{4}\in\N.
\end{equation}

\end{dfn}
In order to have an idea of how frequent they are, we list here all perfect cases for $1\le c\le d\le 10$.
\[\begin{tabular}{ccccc} (1,5) & (1,8) & (2,2) &(2,3) & (2,10) \\ (3,3) &(3,10) & (4,5) & (4,8) &(5,9) \\ (6,6) & (6,7) &(7,7) & (8,9) & (10,10).\end{tabular}\]
If $(c,d)$ is not a perfect case, then $\dim\Sec_r(X)\neq\dim \p^N$ for every $r\in\N$\color{black}. In this case the generic fiber of the map $p_1:\Sec_r(X)\to \p^N$ cannot be zero-dimensional and therefore the general pair of polynomials is not identifiable. Even if $(c,d)$ is a perfect case, the general fiber of $p_1$ may have positive dimension. If this happens, then the general point of $\p^N$ has no preimages under $p_1$, therefore the general pair of ternary forms has no simultaneous decompositions with $r$ summands.

There is another family for which we can easily disprove identifiability \Fra{by looking at the Waring rank of the higher-degree polynomial of the pair. 
}

\begin{pro}\label{pro: se d >> c, allora X è difettiva}
Let $(c,d)$ be a perfect case and let $k$ be the number defined in \eqref{equation:definition of k}. If $\left\lceil\frac13\binom{d+2}{2}\right\rceil>k$, then $X$ is $k$-defective. In particular, the general pair of ternary forms of degrees $c$ and $d$ is not identifiable.
\begin{proof}
Let $(f,g)\in\p^N$ be a general pair of ternary forms of degrees $c$ and $d$. Without loss of generality we assume that $d\neq 2$ and $d\neq 4$. Indeed, there are no perfect cases satisfying our hypothesis for these values of $d$. Since the values $d=2$ and $d=4$ are the only exceptions of 
\Fra{Theorem \ref{thm: AH}}, the number $\left\lceil\frac13\binom{d+2}{2}\right\rceil$ is the Waring rank of $g$. As every simultaneous decomposition of $(f,g)$ gives a decomposition of $g$, the simultaneous rank of $(f,g)$ is at least the rank of $g$. By hypothesis, the latter is strictly greater than $k$. This means that the general pair has no decomposition with $k$ summands. The image of the secant map $p_1:\Sec_k(X)\to\p^N$ has dimension smaller than $N$, hence $X$ is $k$-defective.
\end{proof}
\end{pro}

\begin{rmk}\label{rmk: possiamo assumere che k sia intero} \label{rmk: assumiamo di stare sotto l'iperbole}
In light of our observations, from now on we always suppose that $\binom{c+2}{2}+\binom{d+2}{2}$ is a multiple of 4 and that $k$ is the natural number defined in \eqref{equation:definition of k}. Without loss of generality, we also assume that $p_1:\Sec_k(X)\to \p^N$ is \textit{dominant} - namely, that its image has dimension $N$. 
\Fra{Under these assumptions, the domain and the image of $p_1$ have the same dimension, hence the general fiber of $p_1$ has dimension 0 - we say that $p_1$ is \textit{generically finite}. In other words, $X$ is not $k$-defective and $\SSec_k(X)=\p^N$.} Thanks to Proposition \ref{pro: se d >> c, allora X è difettiva}, this tells us we can work under the assumption that $\left\lceil\frac13\binom{d+2}{2}\right\rceil\le k$. This implies that $\binom{d+2}{2}\le 3k$, that is
\[
	d^2+3d\le 3c^2+9c+4.\]
Under this hypothesis, the only perfect cases with $c\le 2$ are $(2,2)$ and $(2,3)$, the special cases appearing in the statement of Theorem \ref{thm: main}. As we have already observed, they are known to be identifiable since the late nineteenth century. 
For this reason, from now on we suppose that $c\ge 3$. Thanks to  \cite[Theorem 5.1]{AGMO}, we further suppose that $d\neq c+1$.  Under these assumptions, the only cases left in the range $1\le c\le d\le 10$ are \[\begin{tabular}{ccccc} (4,8) &(5,9) &(6,6) &(7,7) &(10,10).\end{tabular}\]
Notice that we excluded the perfect case $(3,3)$. The case of two plane cubics has been classically studied in \cite{London} and \cite{Scorza}. If $c=d=3$, then $k=5$, $N=19$ and $X$ is 5-defective - see also \cite[Remark 4.2]{CarChi03}. In other words, $\dim\SSec_5(X)<19$, thus the secant map $p_1:\Sec_5(X)\to\p^{19}$ cannot be birational.
\end{rmk}
\color{black}

We deal with the map $p_1:\Sec_k(X)\to \p^N$, dominant and generically finite. Our goal is to prove that $p_1$ is not birational, that is, $\deg(p_1)\ge 2$. The following result, proven in \cite[Theorem 5.2]{AGMO}, allows us to reduce the problem.

\begin{thm}\label{thm: tangential projection} Let $V\subset\p^N$ be a nondegenerate irreducible variety of dimension $n$  and let $r\in\N$. 
	Assume that the secant map $p_1:\Sec_{ r}(V)\to\p^N$ is dominant and generically finite. Let $z\in\SSec_{r-1}(V)$ be a general point. Consider the projection $\f:\p^N\dasharrow\p^n$ from the embedded tangent space $\T_z\SSec_{r-1}(V)$. Then
	$\f_{|V}:V\dasharrow\p^n$ is dominant and generically finite, and $\deg(\f_{|V})\le\deg(p_1)$.
\end{thm}
We will apply Theorem \ref{thm: tangential projection} in the case when $V$ is the projective bundle $X$  and $r$ is the number $k$ defined in \eqref{equation:definition of k}. In order to prove that $\deg(p_1)\ge 2$ it is enough to prove that $\deg(\f_{|X})\ge 2$. We want to understand the linear system associated with $\f_{|X}$. 
\Fra{The map} $\f$ described by Theorem \ref{thm: tangential projection} is the projection from the linear space $\T_z\SSec_{k-1}(X)$\Fra{. Once again we apply Lemma \ref{lem: Terracini} to deduce that the linear system associated with $\f_{|X}$ is 
\begin{align*}
    \{&
H\cap X\mid H\subset\p^N
\mbox{ is a hyperplane and }H\supset\T_z\SSec_{k-1}(X)
\}\\
=\{&
H\cap X\mid H\subset\p^N
\mbox{ is a hyperplane and }H\supset\T_{x_i}(X)\mbox{ for every } i\in\{1,\dots,k-1\}\}.
\end{align*}
Keeping in mind that} $H\supset\T_{x_i}(X)$ if and only if $H\cap X$ is singular at $x_i$\Fra{, we are ready to define the linear system we are interested in.}
\begin{dfn}\label{def: sistema LL} Let $X=\p(\oo_{\p^2}(c)\oplus\oo_{\p^2}(d))\subset\p^N$. Let $\Sigma$ be a 0-dimensional subscheme of $X$ consisting of $k-1$ points of multiplicity 2 in general position. We denote by 
$$
\LL:=\p(H^0\oo_X(T_X)\cap I_\Sigma)
$$
the linear system of tautological divisors containing the subscheme $\Sigma$.
\end{dfn}
The linear system $\LL$ on $X$ induces the rational map $\f_{|X}$. As we stressed in Remark \ref{rmk: possiamo assumere che k sia intero}, we work under the assumption that $X$ is not $k$-defective, hence the map $p_1:\Sec_k(X)\to\p^N$ is dominant and generically finite. Thanks to Theorem \ref{thm: tangential projection}, the map $\f_{|X}:X\dashrightarrow\p^3$ is also dominant and generically finite, so the dimension of $\LL$ is 3. Our task is to bound the degree of $\f_{|X}$. 

A standard approach to work with linear systems is to degenerate them. In our case, we will pick some of the points of $\Sigma$ in special position, rather than working with general points. When we perform this kind of degenerations, we have to make sure that we can control the degree of the associated map.
\begin{lem}\label{lem: possiamo specializzare}
	Let $X=\p(\oo_{\p^2}(c)\oplus\oo_{\p^2}(d))\subset\p^N$ and fix $x_1,\ldots,x_{k-1}\in X$. Let $\tilde{\Sigma}$ be the 0-dimensional subscheme of $X$ consisting of $k-1$ points of multiplicity 2 supported at $x_1,\dots,x_{k-1}$. Let $\tilde{\LL}=\p(H^0\oo_X(T_X)\cap I_{\tilde{\Sigma}})$ and call $\tilde{\f}_{|X}$ the associated rational map. If $\dim\tilde{\LL}=\dim\LL=\dim X$, then $\deg(\tilde{\f}_{|X})\le\deg(\f_{|X})$.
\end{lem}
This follows from the more general \cite[Lemma 5.4]{AGMO} and guarantees that we are allowed to degenerate some of the points of $\Sigma$ in special position, as long as our degeneration does not change the dimension of the linear system. In our case, some of the points of $\Sigma$ will belong to a given surface $Z\subset X$. In order to pick a suitable $Z$, consider the bundle morphism
$$
\pi: X\to\p^2
$$ 
and recall that the Picard group of $X$ has rank 2. We choose as generators the tautological divisor $T_X$ and the divisor $\pi^\star (h)$, where $h\subset\p^2$ is a line. We set
$$Z=\{[0,\ell^d]\mid\ell\in\C[x_0,x_1,x_2]_1\}\subset X.$$
It is the section of $\pi$ corresponding to the quotient
\begin{equation}\label{sezione}
\oo_{\p^2}(c)\oplus\oo_{\p^2}(d) \to \oo_{\p^2}(d) \to 0,
\end{equation}
see \cite[Exercise II.7.8]{hartshorne}. In particular, $Z$ is smooth and irreducible and the restriction $\pi _{|Z} : Z \to \p^2$ is an isomorphism. The tautological linear system on $X$ embeds $Z$ as a $d$-Veronese surface. The class of $Z$ is
\begin{equation}\label{classe_zeta}
Z\sim T_X-c\pi^\star (h),
\end{equation}
see for instance \cite[Proposition 9.13]{3264}\footnote{The statement of \cite[Proposition 9.13]{3264} is about a subbundle, instead of a quotient bundle. The reason is that \cite{hartshorne} and \cite{3264} use different conventions, as explained on \cite[page 324]{3264}.}. 

In order to bound $\deg(\f_{|X})$, we want to restrict the map to a suitable subvariety. We will show that this restriction does not increase the degree, provided that such a subvariety is not contained in the contracted locus of $\f_{|X}$.

\begin{dfn}
The \emph{contracted locus} of a rational map is the union of all positive-dimensional fibers. We denote by $\Delta\subset X$ the contracted locus of $\f_{|X}$.\end{dfn}

\begin{lem}
	\label{lem: restrizione a Z non aumenta il grado} Let $f:V\dashrightarrow W$ be a rational map between smooth irreducible varieties. Assume that $f$ is dominant and generically finite. Let $S$ be a subvariety of $V$. If $S$ is not contained in the contracted nor in the indeterminacy locus of $f$, then $\deg(f_{|S})\le\deg(f)$.
\begin{proof}
Let $A=\{p\in W\mid f^{-1}(p)\mbox{ is finite}\}$. By hypothesis $A$ is a nonempty open subset of $W$. Let $B\subset V$ be the indeterminacy locus of $f$ and let $U=f^{-1}(A)\setminus B$. Then $U$ is a nonempty  open subset of $V$. By construction, the restriction $f_{U}:U\to f(U)$ is a finite surjective morphism and $\deg(f_{U})=\deg(f)$. By \cite[Exercise III.9.3(a)]{hartshorne}, $f_{U}$ is flat. By hypothesis $S\cap U\neq\vu$, so $S\cap U$ is a dense open subset of $S$. The flatness of $f_U$ implies that $\deg(f_{U|S})\le\deg(f_U)$, so
\begin{align*}
\deg(f_{|S}) &\le\deg(f_{U|S})\le\deg(f_U)=\deg(f).\qedhere
\end{align*}
\end{proof}
\end{lem}
Now we know that we can bound the degree by specializing some of the base points to $Z$ and then restricting the map. 
Before we proceed, we have to understand a bit better the base points of our linear system $\LL$. Although \color{black} $\Sigma$ is zero-dimensional, the base locus of $\LL$ contains many curves.

\begin{lem}\label{lem: if singular point, then contains the whole fiber}
Let $x\in X$ and let $D\subset X$ be a divisor. If $D\sim T_X$ and $D$ is singular at $x$, then $D\supset \pi^\star (\pi(x))$.
\begin{proof}
Since the class $T_X$ is unisecant, $D\cdot \pi^\star (\pi(x))=1$, so the only possibility for $\mul_x D$ to be greater than 1 is that $D$ contains $\pi^\star (\pi(x))$.
\end{proof}\end{lem}

Concerning our degeneration approach, in order to work, we need to choose carefully the number of points we are going to degenerate. For this reason, we need to check a simple arithmetic property
.

\begin{lem}
	\label{lem: condizione aritmetica} Let $c,d\in\N$. If $\binom{c+2}{2}+\binom{d+2}{2}$ is a multiple of 4, then
\[
3d^2+9d-c^2-3c-12\equiv 0\mod 8.
\]
\begin{proof}
By hypothesis there exists $t\in\N$ such that $\binom{c+2}{2}+\binom{d+2}{2}=4t$. This means that
\begin{align*}
c^2&+3c+d^2+3d+4=8t\Rightarrow 3d^2+9d+3c^2+9c+12=24t,
\end{align*}
hence $3d^2+9d-c^2-3c-12
=8(3t-3)-4c(c+3)\equiv 0\mod 8$.\qedhere
\end{proof}
\end{lem}
We are actually interested in the class of $3d^2+9d-c^2-3c-12$ modulo 16. By Lemma \ref{lem: condizione aritmetica}, there are two possibilities: either $3d^2+9d-c^2-3c-12\equiv 0\mod 16$ or $3d^2+9d-c^2-3c-12\equiv 8\mod 16$. We use two different strategies to deal with these two cases.

\section{The first case}\label{section: first case}
The goal of this section is to prove Theorem \ref{thm: main} when
\begin{equation}\label{eq: primo caso modulo 16}
3d^2+9d-c^2-3c-12\equiv 0\mod 16.
\end{equation}
This is accomplished in  Proposition \ref{pro: dimostrazione nel primo caso}. We start by setting up the degeneration we need. Under hypothesis \eqref{eq: primo caso modulo 16} we define the integer
\[s_1=\frac{3d^2+9d-c^2-3c-12}{16}.\]
We want to make sure that \Fra{$s_1\in\{0,\ldots,k-1\}$}. Thanks to Remark \ref{rmk: assumiamo di stare sotto l'iperbole} we work under the assumption $d\ge c\ge 3$, so
\begin{align*}
s_1\ge\frac{2c^2+6c-12}{16}\ge 0.
\end{align*}
Remark \ref{rmk: assumiamo di stare sotto l'iperbole} also allows us to assume that $d^2+3d\le 3c^2+9c+4$\color{black}, hence
\begin{align*}
k-1-s_1&=\frac{c^2+3c+d^2+3d-4}{8}-\frac{3d^2+9d-c^2-3c-12}{16}
&=\frac{3c^2+9c-d^2-3d+4}{16}\ge 0\color{black}.
\end{align*}
Consider the linear system $\LL$ on $X$ introduced in Definition \ref{def: sistema LL}. Since \Fra{$s_1\in\{0,\ldots,k-1\}$}, we can degenerate $s_1$ of the $k-1$ base points of $\LL$ in special position.
\begin{dfn}\label{def: specializzazione, primo caso}
Let $X=\p(\oo_{\p^2}(c)\oplus\oo_{\p^2}(d))\subset\p^N$ and let $Z\subset X$ be the section defined by \eqref{sezione}. Let $z_1,\dots,z_{s_1}$ be general points of $Z$ and let $x_{s_1+1},\dots,x_{k-1}$ be general points of $X$; in particular, $x_{s_1+1},\dots,x_{k-1} \not\in Z$. Let $\Sigma_1$ be the 0-dimensional subscheme of $X$ consisting of $k-1$ points of multiplicity 2 supported at $z_1,\dots,z_{s_1},x_{s_1+1},\dots,x_{k-1}$. Similarly to Definition \ref{def: sistema LL}, we define $
L_1:=H^0\oo_X(T_X)\cap I_{\Sigma_1}
$ and $\LL_1:=\p(L_1)$. We call $\f_1$ the associated rational map.
\end{dfn}
In order to be able to perform the computations on $\LL_1$, we need to prove that this specialization does not increase the dimension of the linear system. In our case, the Castelnuovo exact sequence \eqref{exact_sequence} becomes
\begin{equation}\label{Castelnuovo, first case}
0\to L_1\cap I_Z\to L_1\to(L_1)_{Z}.
\end{equation}
Now we are in position to describe both the right and the left hand sides of sequence \eqref{Castelnuovo, first case} and to find their dimensions, thereby computing $\dim(\LL_1)$.
\begin{lem}\label{lem: primo caso, la dimensione non cresce}
In the specialization of Definition \ref{def: specializzazione, primo caso}, we have
\begin{enumerate}
\item \label{item: primo caso, solo una sezione contenente Z}$  L_1\cap I_Z\cong L_2^c(2^{k-1-s_1},1^{s_1})$ and it has dimension 1,
\item\label{item: primo caso, la restrizione ha dim=3} 
$( L_1)_Z\cong L_2^d(2^{s_1},1^{k-1-s_1})$
and it has dimension 3,
\item\label{item: L1 mappa su P3} $\dim( L_1)=4$,
\item\label{item: L1 ristretto è quello giusto} $( L_1)_Z=L_{1|Z}$.
\end{enumerate}
\begin{proof}We will use the bundle morphism $\pi:X\to\p^2$ to translate the question on our linear systems in terms of linear systems on the plane. Recall that $\pi$ restricts to an isomorphism between $Z$ and $\p^2$. 
Let us start by showing that $L_1\cap I_Z$ is isomorphic to a vector subspace of $L_2^c(2^{k-1-s_1},1^{s_1})$ of dimension at most 1. Take $D\in\LL_1$ containing $Z$. Then there exists a divisor $D^\prime$ such that $D=Z+D^\prime$. From \eqref{classe_zeta} we obtain $D^\prime\sim c\pi^\star (h)$,
so it projects to a plane curve of degree $c$. Notice that the general $D$ has multiplicity 2 at $z_1,\dots,z_{s_1}, x_{s_1+1}
,\dots, x_{k-1}$, while $Z$ has multiplicity 1 at $z_1,\dots,z_{s_1}$ and does not contain $x_{s_1+1}
,\dots, x_{k-1}$. It follows that $D^\prime$ has multiplicity 1 at $z_1,\dots,z_{s_1}$ and multiplicity 2 at $x_{s_1+1}
,\dots x_{k-1}$. The correspondence $D\mapsto D^\prime$ is an isomorphism, therefore, after the projection on $\p^2$, we can regard elements of $\p(L_1\cap I_Z)$ as plane curves of degree $c$ singular at $\pi(x_{s_1+1})
,\dots, \pi(x_{k-1})$ and passing through $\pi(z_1),\dots,\pi(z_{s_1})$. Hence $L_1\cap I_Z$ is 
%
%
a vector subspace of $L_2^c(2^{s_1},1^{k-1-s_1})$. The latter has dimension 1 by Theorem \ref{thm: AH}, 
hence $\dim(L_1\cap I_Z)\le 1$.

In a similar way, now we prove that $( L_1)_Z$ is a vector subspace of $L_2^d(2^{s_1},1^{k-1-s_1})$ of dimension at most 3. Elements of $(\LL_1)_Z$ have class $T_{X|Z}$, are singular at $s_1$ general points and pass through $k-1-s_1$ simple base points in general position. Indeed, by Lemma \ref{lem: if singular point, then contains the whole fiber} the base locus $\Bs(\LL_1)$ contains not only $\Sigma_1$, but also $k-1$ general fibers. Since $Z$ is a section of $\pi$, each of these fibers intersects $Z$ in one point. Therefore curves in $(\LL_1)_{Z}:=\p((L_1)_Z)$ are not only singular at $z_1,\dots,z_{s_1}$, but they also contain $k-1-s_1$ simple base points in general position, namely the intersections of $Z$ with the fibers $\pi^{-1}(\pi(x_{s_1+1})),\dots,\pi^{-1}(\pi(x_{k-1}))$. 

Now we show that the linear system $|T_{X|Z}|$
corresponds isomorphically to $|dh|$ via the morphism $\pi$. Denote by $c_1(\oo_{\p^2}(c)\oplus\oo_{\p^2}(d))$ and $c_2(\oo_{\p^2}(c)\oplus\oo_{\p^2}(d))$ the first and the second Chern classes of the line bundle $\oo_{\p^2}(c)\oplus\oo_{\p^2}(d)$ - a good reference on Chern classes is \cite[Chapter 5]{3264}. Recall that by the Whitney formula we have
$$
c_1(\oo_{\p^2}(c)\oplus\oo_{\p^2}(d)) = (c+d) h\mbox{ and }
c_2(\oo_{\p^2}(c)\oplus\oo_{\p^2}(d))=c d(\mbox{pt}),
$$
where $\mbox{pt}$ indicates the class of a point of $\pp^2$.
The fundamental relation
\begin{align*}
T_X^2&\equiv T_X\cdot\pi^\star (c_1(\oo_{\p^2}(c)\oplus\oo_{\p^2}(d)))-\pi^\star (c_2(\oo_{\p^2}(c)\oplus\oo_{\p^2}(d)))\\
&=(c+d)T_X\cdot\pi^\star (h)-cd\pi^\star (\mbox{pt}),
\end{align*}
discussed for instance in \cite[Section A.3]{hartshorne}, implies that
$$
T_X \cdot Z=T_X^2 - c \pi^\star (h) \cdot T_X=d \pi^\star (h) \cdot T_X - cd\pi^\star (\mbox{pt}).
$$
It follows that $(\pi_{|Z})_\star T_X \cdot Z \sim dh$, therefore $(L_1)_Z$ is isomorphic to a vector subspace of $L_2^d(2^{s_1},1^{k-1-s_1})$. Again by 
\Fra{Theorem \ref{thm: AH}}, the latter has dimension 3, so $\dim((L_1)_Z)\le 3$. Since $\dim(L_1\cap I_Z)\le 1$ and $\dim((L_1)_Z)\le 3$, exact sequence \eqref{Castelnuovo, first case} gives $$\dim L_1\le\dim( L_1\cap I_Z)+\dim(( L_1)_Z)\le 1+3=4.$$ On the other hand, $\LL_1$ is a degeneration of $\LL=\p(L)$, so $\dim L_1\ge\dim L\ge 4$ by semicontinuity. This proves the third claim.

As a consequence, $ L_1\cap I_Z$ and $( L_1)_Z$ have indeed dimension 1 and 3, respectively, otherwise $\dim( L_1)$ would be smaller than 4. This proves the first two claims.

We are left to prove the fourth part. Since $\dim( L_1)-\dim(L_1\cap I_Z)=\dim(L_1)_Z$, the rightmost arrow in the sequence \eqref{Castelnuovo, first case} is surjective, hence $(L_1)_Z=L_{1|Z}$. \end{proof}
\end{lem}
Now we know that \Fra{the specialization we introduced in Definition \ref{def: specializzazione, primo caso} preserves the dimension of the linear system. In other words, $\dim\LL_1=\dim\LL$, so the codomain of $\f_1$ is $\p^3$ and we can take advantage of Lemma \ref{lem: possiamo specializzare}. The following result proves Theorem \ref{thm: main} under hypothesis \eqref{eq: primo caso modulo 16}.} 

\begin{pro}\label{pro: dimostrazione nel primo caso}The 
map \Fra{$\f_1:X\dashrightarrow\p^3$ }  associated with $\LL_1$ is not birational.
	\begin{proof}
We want to apply Lemma \ref{lem: restrizione a Z non aumenta il grado} and show that  $\f_{1|Z}:Z\dashrightarrow\p^2$ is not birational. For this purpose, the first thing we need is to show that $Z\not\subset\Delta$. Assume by contradiction that $Z\subset\Delta$. By Lemma \ref{lem: primo caso, la dimensione non cresce}\eqref{item: primo caso, la restrizione ha dim=3}, the image of $Z$ is a nondegenerate plane curve $Y\subset\p^2$, and images of divisors in $\LL_{1|Z}$ are line sections of $Y$. This implies that the general element of $\LL_{1|Z}$ is reducible. It follows by Lemma \ref{lem: primo caso, la dimensione non cresce}\eqref{item: L1 ristretto è quello giusto} that   the general element of $\LL_2^d(2^{s_1},1^{k-1-s_1})$ is reducible, and this contradicts Lemma \ref{lemma: arbarellocornalba}. The case $(6,9,0)$ cannot happen because \Fra{$(2,6)$ is not perfect.} 
Moreover, since $\dim L_1 =4$ by Lemma \ref{lem: primo caso, la dimensione non cresce}\eqref{item: L1 mappa su P3}, we see that $Z$ is not in the base locus of $\LL_1$ by Lemma \ref{lem: primo caso, la dimensione non cresce}\eqref{item: primo caso, solo una sezione contenente Z}, hence it is not contained in the indeterminacy locus of $\f_1$.

Now we only have to prove that $\f_{1|Z}:Z\dashrightarrow\p^2$ is not birational. It suffices to show that the general element of $\LL_2^d(2^{s_1},1^{k-1-s_1})$ is not a rational curve. Thanks to Lemma \ref{lemma: arbarellocornalba}, the general element of $\LL_2^d(2^{s_1},1^{k-1-s_1})$ is irreducible 
and is singular at exactly $s_1$ ordinary double points, so it has genus
\begin{align*}
    \frac{(d-1)(d-2)}{2}-s_1=\frac{5d^2-33d+c^2+3c+28}{16}.
\end{align*}
Under the assumptions we made in Remarks \ref{rmk: possiamo assumere che k sia intero}, together with hypothesis \eqref{eq: primo caso modulo 16}, it is not restrictive to suppose that $c\ge 4$ and $d\ge 6$. Then we can bound the genus as
\[\frac{5d^2-33d+c^2+3c+28}{16}\ge\frac{5d^2-33d+56}{16}>0.\qedhere\]
	\end{proof}
\end{pro}

\section{The second case}\label{sec:secondo caso}
The goal of this section is to prove Theorem \ref{thm: main} when
\begin{equation}\label{eq: secondo caso modulo 16}
3d^2+9d-c^2-3c-12\equiv 8\mod 16.
\end{equation}
This is accomplished in Proposition \ref{pro: dimostrazione nel secondo caso}. Under the assumptions we made in Remark \ref{rmk: possiamo assumere che k sia intero}
, the only case satisfying \eqref{eq: secondo caso modulo 16} with $c\le 7$ is $(3,3)$ which, as we have already observed, cannot be identifiable. \color{black}  Therefore in this section we assume that $c\ge 8$. As in Section \ref{section: first case}, we start by setting up the degeneration we need. Under hypothesis \eqref{eq: secondo caso modulo 16} we can define the integer
\[s_2=\frac{3d^2+9d-c^2-3c-4}{16}.\]
Just as we did in Section \ref{section: first case} for $s_1$, it is easy to check that $s_2\in\{0,\dots,k-1\}$.
Consider the linear system $\LL$ introduced in Definition \ref{def: sistema LL}. 
\Fra{We degenerate} $s_2$ of the $k-1$ base points of $\LL$ in special position.
\begin{dfn}\label{def: specializzazione, secondo caso}
Let $X=\p(\oo_{\p^2}(c)\oplus\oo_{\p^2}(d))\subset\p^N$ and let $Z\subset X$ be the section given by \eqref{sezione}. Let $z_1,\dots,z_{s_2}$ be general points of $Z$ and let $x_{s_2+1},\dots,x_{k-1}$ be general points of $X$. Let $\Sigma_2$ be the 0-dimensional subscheme of $X$ consisting of $k-1$ points of multiplicity 2 supported at $z_1,\dots,z_{s_2},x_{s_2+1},\dots,x_{k-1}$. We define $L_2:=H^0\oo_X(T_X)\cap I_{\Sigma_2}$ and $\LL_2=\p(L_2)$. We call $\f_2$ the associated rational map.
\end{dfn}
Again, our first concern is to check that the degeneration presented in Definition \ref{def: specializzazione, secondo caso} does not increase the dimension of the linear system. In other words, we want to prove that $\dim( L_2)=4$. 
\begin{lem}\label{lem: secondo caso, la dimensione non cresce}
In the specialization of Definition \ref{def: specializzazione, secondo caso}, we have
	\begin{enumerate}
\item\label{item: secondo caso, la restrizione ha dim=2} $( L_2)_Z\cong L_2^d(2^{s_2},1^{k-1-s_2})$ and it has dimension 2,
\item \label{item: secondo caso, residuo è un pencil}$ L_2\cap I_Z\cong L_2^c(2^{k-1-s_2},1^{s_2})$ and it has dimension $2$,
\item\label{item: L2 mappa su P3} $\dim( L_2)=4$,
\item\label{item: L2 ristretto è quello giusto} $( L_2)_Z= L_{2|Z}$.
	\end{enumerate}
	\begin{proof}
The proof goes exactly as in Lemma \ref{lem: primo caso, la dimensione non cresce}. We only have to check that both \Fra{$L_2^d(2^{s_2},1^{k-1-s_2})$ and $L_2^c(2^{k-1-s_2},1^{s_2})$ have dimension} 2.
	\end{proof}
\end{lem}

\begin{rmk} As a byproduct of Lemmas \ref{lem: primo caso, la dimensione non cresce} and \ref{lem: secondo caso, la dimensione non cresce}, we obtain that $\dim\LL=3$, even without assuming so from the beginning. This means that, whenever $(c,d)$ is a perfect case, $k-1$ general double points impose independent conditions on $\oo_X(T_X)$. By Lemma \ref{lem: Terracini} \Fra{and equation \eqref{eq: sistemi lineari vs secanti}}, $X$ is not $(k-1)$-defective.\end{rmk}

The main difference between situation \eqref{eq: primo caso modulo 16} and situation \eqref{eq: secondo caso modulo 16} is that in this second case $(\LL_2)_Z$ induces a map to $\p^1$, instead of $\p^2$. Therefore $Z\subset \Delta$, in this setting. This means that we cannot apply Lemma \ref{lem: restrizione a Z non aumenta il grado} to $Z$, but rather we will find another suitable subvariety. Observe that if $T\in\p(L_2\cap I_Z)$, then $T\sim Z+V$ for some element $V\sim c\pi^*(h)$. By Lemma \ref{lem: secondo caso, la dimensione non cresce}\eqref{item: secondo caso, residuo è un pencil}, there is a pencil of such $V$'s.

\begin{lem}
\label{lem: buone proprietà di T e V}
Let $B$ be a general element of $\LL_2^c(2^{k-1-s_2},1^{s_2})$ and $V=\pi^\star B\subset X$. Let $T$ be a general element of $\LL_2$. Then
\begin{enumerate}
\item \label{item: V irriducibile e mappa su P2}$V$ is irreducible, it is not contained in the contracted locus $\Delta$ and $\overline{\f_2(V)}=\p^2$.
\item $T\not\supset V$ and $T\cap V\not\subset\Delta$. Moreover $\overline{\f_2(T\cap V)}=\p^1$.
\end{enumerate}
\begin{proof}
Since $V$ corresponds to a general element of $\LL_2^c(2^{k-1-s_2},1^{s_2})$, it is irreducible by Lemma \ref{lemma: arbarellocornalba}. By construction we have $Z\subset\Delta$. As $\codim \Delta \ge 1$, we can choose $V$ so that $V \not \subset \Delta$. Moreover, we can also assume that $V$ is not contained in the indeterminacy locus of $\f_2$.
Hence $\dim(\f_2(V))=2$. Observe that $V+Z\in\LL_2$, so $\overline{\f_2(V\cup Z)}=\p^2$. Since $\dim(\f_1(Z))=1$, the image of $V$ is $\p^2$. This completes the proof of the first claim.

Now assume by contradiction that the general element $T$ of $\LL_2$ contains $V$. Then $L_2\cap I_Z$ would have only one element, up to scalar. Namely, $L_2\cap I_Z=\langle V + Z\rangle$. This contradicts Lemma \ref{lem: secondo caso, la dimensione non cresce}\eqref{item: secondo caso, residuo è un pencil}. Since $T$ does not contain $V$, the intersection is a curve on $V$. If such a curve was contained in $\Delta$ for a general $T\in\LL_2$, then $V\subset\Delta$, in contradiction to part \eqref{item: V irriducibile e mappa su P2}. Finally, $T\cap V$ is an element of $\LL_{2|V}$. AS $\overline{\f_2(V)}=\p^2$, we have $\overline{\f_2(T\cap V)}=\p^1$.
\end{proof}
\end{lem}


We are interested in $T\cap V$. This is not an irreducible curve, because both $T$ and $V$ contain the $k-1$ fibers of the base locus of $\LL_2^c(2^{k-1-s_2},1^{s_2})$ via $\pi$. However, $T\sim T_X$ is a unisecant divisor, that is $T$ either contains a fiber of $\pi$ or $T$ intersects 
a fiber precisely in one point, hence there exists a unique horizontal component of $T\cap V$. We define $C$ to be such an irreducible component. As $T$ is general, $C\not\subset\Delta$ and $C$ is not contained in the indeterminacy locus of $\f_2$. Our strategy to bound $\deg(\f_2)$ is to restrict the map to $C$.
\begin{pro}
    \label{pro: dimostrazione nel secondo caso}
The 
map \Fra{$\f_2:X\dashrightarrow\p^3$ } associated with $\LL_2$ is not birational.
\begin{proof}
%
As in Lemma \ref{lem: buone proprietà di T e V}, we define $B$ as a general element of $\LL_2^c(2^{k-1-s_2},1^{s_2})$ and $T$ a general element of $\LL_2$. We want to apply Lemma \ref{lem: restrizione a Z non aumenta il grado} and show that $\f_{2|C}:C\dashrightarrow\p^1$ is not birational. We
only need to prove that $C$ is not a rational curve. Since $\pi_{|T}$ is birational, it restricts to a birational map between $C$ and $B$, hence  they have the same genus. By Lemma \ref{lemma: arbarellocornalba}, the curve $B$ has only ordinary double points, so its genus is
\begin{align*}
\binom{c-1}{2}-(k-1-s_2)
=\frac{5c^2+d^2-33c+3d+20}{16}\ge \frac{6c^2-30c+20}{16}\ge 1
\end{align*}
for every $c\ge 8$.
\end{proof}
\end{pro}

\bibliographystyle{alpha}
\bibliography{simultaneousreferences.bib}

\begin{thebibliography}{AGMO18}

\bibitem[AC81]{AC}
E.~Arbarello and M.~Cornalba.
\newblock Footnotes to a paper of {B}eniamino {S}egre: ``{O}n the modules of
  polygonal curves and on a complement to the {R}iemann existence theorem''
  ({I}talian) [{M}ath. {A}nn. {\bf 100} (1928), 537--551; {J}buch {\bf 54},
  685].
\newblock {\em Mathematische Annalen}, 256(3):341--362, 1981.

\bibitem[AGMO18]{AGMO}
E.~Angelini, F.~Galuppi, M.~Mella, and G.~Ottaviani.
\newblock On the number of {W}aring decompositions for a generic polynomial
  vector.
\newblock {\em Journal of Pure and Applied Algebra}, 222(4):950--965, 2018.

\bibitem[BBC12]{BalBerCat12}
E.~Ballico, A.~Bernardi, and M.~V. Catalisano.
\newblock Higher secant varieties of $\mathbb{P}^n\times\mathbb{P}^1$ embedded
  in bi-degree $(a,b)$.
\newblock {\em Communications in algebra}, 40(10):3822--3840, 2012.

\bibitem[Cam92]{Campbell}
J.~E. Campbell.
\newblock Note on the maximum number of arbitrary points which can be double
  points on a curve, or surface, of any degree.
\newblock {\em The Messenger of Mathematics}, XXI:158--164, 1891-1892.

\bibitem[CC02]{trisecant}
L.~Chiantini and C.~Ciliberto.
\newblock Weakly defective varieties.
\newblock {\em Transactions of the American Mathematical Society},
  354(1):151--178, 2002.

\bibitem[CC03]{CarChi03}
E.~Carlini and J.~Chipalkatti.
\newblock On {W}aring's problem for several algebraic forms.
\newblock {\em Commentarii Mathematici Helvetici}, 78(3):494--517, 2003.

\bibitem[CR06]{CR}
C.~Ciliberto and F.~Russo.
\newblock Varieties with minimal secant degree and linear systems of maximal
  dimension on surfaces.
\newblock {\em Advances in Mathematics}, 200(1):1--50, 2006.

\bibitem[Dal84]{dale}
M.~Dale.
\newblock Terracini's lemma and the secant variety of a curve.
\newblock {\em Proceedings of the London Mathematical Society. Third Series},
  49(2):329--339, 1984.

\bibitem[DDL14]{REF}
I.~Domanov and L.~De~Lathauwer.
\newblock Canonical polyadic decomposition of third-order tensors: reduction to
  generalized eigenvalue decomposition.
\newblock {\em SIAM Journal on Matrix Analysis and Applications},
  35(2):636--660, 2014.

\bibitem[EH16]{3264}
D.~Eisenbud and J.~Harris.
\newblock {\em 3264 and all that---a second course in algebraic geometry}.
\newblock Cambridge University Press, Cambridge, 2016.

\bibitem[GM19]{GM}
F.~Galuppi and M.~Mella.
\newblock Identifiability of homogeneous polynomials and {C}remona
  transformations.
\newblock {\em Journal f\"{u}r die reine und angewandte Mathematik},
  757:279--308, 2019.

\bibitem[Har77]{hartshorne}
R.~Hartshorne.
\newblock {\em Algebraic geometry}.
\newblock Springer-Verlag, New York-Heidelberg, 1977.
\newblock Graduate Texts in Mathematics, No. 52.

\bibitem[KMM21]{MOUR}
R.~Khouja, P.-A. Mattei, and B.~Mourrain.
\newblock Tensor decomposition for learning {G}aussian mixtures from moments,
  2021.
\newblock Preprint arXiv:2106.00555.

\bibitem[Lan12]{Lan}
J.~M. Landsberg.
\newblock {\em Tensors: geometry and applications}, volume 128 of {\em Graduate
  Studies in Mathematics}.
\newblock American Mathematical Society, Providence, RI, 2012.

\bibitem[Lon90]{London}
F.~London.
\newblock {\"{U}}ber die {P}olarfiguren der ebenen {C}urven dritter {O}rdnung.
\newblock {\em Mathematische Annalen}, 36(4):535--584, 1890.

\bibitem[LRA93]{LRA}
S.~E. Leurgans, R.~T. Ross, and R.~B. Abel.
\newblock A decomposition for three-way arrays.
\newblock {\em SIAM Journal on Matrix Analysis and Applications},
  14(4):1064--1083, 1993.

\bibitem[Mel06]{M1}
M.~Mella.
\newblock Singularities of linear systems and the {W}aring problem.
\newblock {\em Transactions of the American Mathematical Society},
  358(12):5523--5538, 2006.

\bibitem[Mel09]{M2}
M.~Mella.
\newblock Base loci of linear systems and the {W}aring problem.
\newblock {\em Proceedings of the American Mathematical Society},
  137(1):91--98, 2009.

\bibitem[OS10]{OS}
G.~Ottaviani and E.~Sernesi.
\newblock On the hypersurface of {L}\"{u}roth quartics.
\newblock {\em Michigan Mathematical Journal}, 59(2):365--394, 2010.

\bibitem[Rob89]{Roberts}
R.~A. Roberts.
\newblock Note on the plane cubic and a conic.
\newblock {\em Proceedings of the London Mathematical Society}, 21, 1889.

\bibitem[Sco99]{Scorza}
G.~Scorza.
\newblock Sopra le figure polari delle curve piane del terzo ordine.
\newblock {\em Mathematische Annalen}, 51:154--157, 1899.

\bibitem[Ter11]{terracini}
A.~Terracini.
\newblock Sulle {$V_k$} per cui la varietà degli {$S_h$} $(h+1)$-seganti ha
  dimensione minore dell'ordinario.
\newblock {\em Rendiconti del Circolo Matematico di Palermo (1884-1940)},
  31:392--396, 1911.

\bibitem[Ter15]{terracinicoppiediternarie}
A.~Terracini.
\newblock Sulla rappresentazione delle coppie di forme ternarie mediante somme
  di potenze di forme lineari.
\newblock {\em Annali di Matematica Pura ed Applicata (1898-1922)}, 24:1--10,
  1915.

\end{thebibliography}
\end{document}